\documentclass{class}

\usepackage{amsmath,amsfonts,amssymb,color,enumerate}

\def\N{\mathbb N}  
  
\def\R{\mathbb R}

\def\Z{\mathbb Z}  
 
\def\phi{\varphi}
\def\.{{\cdot}}  
\def\<{\langle}  
\def\>{\rangle}  
\def\({\big(}   
\def\){\big)}

\def\implies{\hbox{${\Rightarrow}$}}    
\def\<{\langle} 
\def\>{\rangle}
\def\defi{\stackrel{\rm def}{=}}

\def\ssk{\smallskip}  
\def\msk{\medskip}   
\def\bsk{\bigskip}   
\def\nin{\noindent}  
\def\cen{\centerline}

\def\lead{\leaders\hbox to 1.5ex{\hss${.}$\hss}\hfill}
\def\arr{\hbox to 60pt{\rightarrowfill}}
\def\larr{\hbox to 60pt{\leftarrowfill}}

\theoremstyle{defn}
\theoremstyle{thm}

\long\def\alert#1{\parindent2em\smallskip\hbox to\hsize%
{\hskip\parindent\vrule%
\vbox{\advance\hsize-2\parindent\hrule\smallskip\parindent.4\parindent%
\narrower\noindent#1\smallskip\hrule}\vrule\hfill}\smallskip\parindent0pt}

\title{On the Component Factor Group $G/G_0$\\ 
                    of a Pro-Lie Group $G$}
                  
\author{Rafael Dahmen and Karl H. Hofmann}

\lastname{Dahmen and Hofmann}


         
\address{%
Rafael Dahmen\\
Karlsruher Institut f\"ur Technologie\\ 
(KIT)\\
76131 Karlsruhe, Germany\\
rafael.dahmen@kit.edu}

\address{%
Karl Heinrich Hofmann\\       
Fachbereich Mathematik\\      
Technische Universit\"at Darmstadt\\
Schlossgartenstra{\ss}e 7\\
64289 Darmstadt, Germany\\       
hofmann@mathematik.tu-darmstadt.de}

\begin{document}

\maketitle

\begin{abstract}
	A pro-Lie group $G$ is a topological group such that
	$G$ is isomorphic to the projective limit of all 
	quotient groups $G/N$ (modulo closed normal subgroups $N$)
	such that $G/N$ is a finite dimensional real Lie group.
	A topological group is almost connected if the totally
	disconnected factor group $G_t\defi G/G_0$ of $G$ modulo the
	identity component $G_0$ is compact. 
        In this case it is straightforward that each Lie group
        quotient $G/N$ of $G$ has finitely many components. 
        However, in spite of
	a comprehensive literature on pro-Lie groups, the following
	theorem, proved here, was not available until now:
	\hfill\break
	{\sc Theorem}. {\it A pro-Lie group $G$ is almost connected if
		each of its Lie group quotients $G/N$ has finitely 
		many connected components.}
		\hfill\break
	\nin The difficulty of the proof is the verification of
	the completeness of $G_t$.
	\end{abstract}

\bsk

\section*{Projective Limits of Almost Connected Lie Groups}

A notorious problem in the structure theory of  { pro-Lie 
groups} is
the completeness
of quotient groups, notably that of the group $G/G_0$ of connected components.
In one of the sources on pro-Lie groups,
\cite{probook}, the { section} following 
Definition 4.24 on pp.195ff.
exhibits some of the characteristic difficulties involving the
completeness of quotients of a pro-Lie group $G$ in general and the 
quotient $G_t=G/G_0$ in particular. In their
entirety, these difficulties remain unresolved today. We shall settle
the completeness issue of $G_t$ here for any pro-Lie group whose
Lie group quotients have finitely many components.

\msk
{
	 Existing literature (see \cite{almost}, Corollary 8.4)
          provides the  following conclusion, which reinforces the 
         independent interest in the result of this note:

\msk

	{\it An almost connected pro-Lie group $G$ contains a 
         maximal compact subgroup $C$
		and a closed subspace $V$
	         homeomorphic to $\R^J$ for a set $J$
                such  that 

\cen{$(c,x)\mapsto cx\colon C\times V\to G$}
		\vskip-2pt\nin is a homeomorphism.}}

\bsk

So, let $G$  denote a topological group
and ${\cal N}(G)$  the set of all closed normal
subgroups of $G$ for which $G/N$ is a Lie group.
With these conventions we formulate a theorem,
to be proved subsequently.
The proof  requires some effort. It is based
on information from \cite{probook}.

\msk

\begin{Theorem} \label{th:main} For a pro-Lie group $G$, the following
	statements are equivalent:
	\begin{enumerate}[\rm(1)]
	
		\item $G_t$ is compact,
	
		\item There is a compact totally disconnected
 		subspace $D\subseteq G$
		being mapped homeo\-morphically onto $G_t$ by the 
		quotient map $q_t\colon G\to G_t$.
		
		\item $(G/N)_t$ is finite for all $N\in{\cal N}(G)$.
		\end{enumerate}
\end{Theorem}

\nin The proof of the  theorem will require the proof of some new lemmas
and some references to existing literature. 
The first one is { cited from \cite{almost}, 
Main Theorem 8.1, Corollary 8.3.} 

\begin{Lemma} \label{l:one} Let $G$ be an almost connected
	pro-Lie group. Then the following conclusions hold:
	
	\begin{enumerate}[\rm(i)]
	
		\item $G$ contains a maximal compact subgroup $C$, 
		and any compact subgroup of $G$ has a conjugate inside $C$.

		\item $G=G_0C$. 

		\item $G$ contains a profinite subgroup $D$ such that
                $G=G_0D$.
	\end{enumerate}
\end{Lemma}

\msk

\nin
For every compact group $K$ there is a compact
totally disconnected subspace $D\subseteq K$ such that
$(k,d)\mapsto kd:K_0\times D\to K$ is a homeomorphism
(see \cite{book3}, Corollary 10.38, p.~573). From 
Lemma \ref{l:one}
we know that $G$ is almost connected iff 
there is a compact subgroup $K\subseteq G$
such that $G=G_0K$. Write $K=K_0D$ with the topological
direct factor $D\subseteq K$ as we just pointed out. Then
$G=G_0K_0D=G_0D$ and so $(g,d)\mapsto gd:G_0\times D\to G$
is readily seen to be a homeomorphism. Thus, in
Theorem \ref{th:main}, Condition (1) implies (2),
and for (2)\implies(1) there is nothing to prove. 

\msk

Let us establish that (1)\implies (3): 

\nin
Assume $G/N$ to be a Lie group quotient of $G$.
Then $(G/N)_0=G_0N/N\cong G_0/(G_0\cap N)$ (cf.\ \cite{probook},
Lemma 3.29, p.152). Let $K$ be a compact subgroup of $G$ 
such that $G=G_0K$, and let $L=G/N$. Then $L=(G_0N/N)(KN/N)
=L_0C$ for the compact Lie group $C=KN/N$. Thus
$L_t=L_0C/L_0\cong C/(C\cap L_0)$ is a compact totally disconnected
Lie group and is therefore finite. This proves (3).

\bsk

There remains a proof of the implication (3)\implies (1): 

\nin
For the moment let us assume that the
following hypothesis is satisfied

\msk 

\cen{
{\bf (H)} $G_t$ \it is a complete topological group.}

\msk

\nin 
By \cite{probook}, Corollary 3.31, hypothesis (H) implies that
$G_t$ is prodiscrete, that is, 
$G_t=\lim_{N\in{\cal N}(G_t)} G_t/N$ where $G_t/N$ is discrete.
Now $G_t/N$ is a Lie group quotient of $G$ and so is finite
by (3). Hence $G_t$ is profinite and thus compact. This
proves Condition (1).

\msk

It therefore remains to prove (H). 
For this purpose we shall invoke results from \cite{probook}, pp.195ff.

\msk

Firstly, we define ${\cal M}(G)$ to be the subset of
all $M\in{\cal N}(G)$ with the additional property that 
each open subgroup $N\subseteq M$ from ${\cal N}(G)$ has
finite index in $M$. We shall then use

\begin{Lemma} \label{l:two} If $G$ is a pro-Lie group such that

	\nin{\bf(*)} each Lie group quotient $G/N$, ${ N}
        \in{\cal N}(G)$ is almost connected, 

	\nin
	then ${\cal M}(G)$ is cofinal in 
        ${\cal N}({G})$ and thus is a filter basis. 
	
	Moreover, $G$ is the strict projective limit 
	of the $G/M$, $M\in{\cal M}(G)$.
\end{Lemma}

\begin{Proof} 
	For the proof see Lemma
	4.25 in \cite{probook}, pp.195 and 196.
\end{Proof}

\msk

\nin We note that Lemma 4.25 in \cite{probook} states as hypothesis
that $G$ is almost connected which implies {\bf(*)}. But 
the hypothesis {\bf(*)} is all that is used in the proof
of Lemma 4.25.

\bsk 

\nin
Any set $\cal Z$ of subsets of a set 
$G$ may be considered as a subbasis { of closed sets for}
a topology. If $G$ is { a topological group}, and 
$\cal Z$ is the set of all cosets $gM$ with $g\in G$ and 
$M\in{\cal M}(G)$,
then $\cal Z$ generates { the set of closed sets of}
a topology on $G${,} called the Z-{\it topology}.

\msk

\begin{Lemma} \label{l:three} The {\rm Z}-topology on a
	pro-Lie group $G$ satisfying {\rm Condition} {\bf(*)}
	of {\rm Lemma \ref{l:two}} is a compact $T_1$-topology.
\end{Lemma}

\begin{Proof} See Proposition 4.27 of \cite{probook}, pp.~197--201.
\end{Proof}

\msk

\nin Again we note that Proposition 4.27 { in \cite{probook}}
  assumes the hypothesis that
$G$ is almost connected, but the proof of the conclusion of
Lemma \ref{l:three} only uses Hypothesis {\bf(*)} of 
Lemma \ref{l:two}.

\msk

\nin
We now adjust the proof of Theorem 4.28 on p.~202 of \cite{probook}
for our purposes.

\msk

\begin{Lemma} \label{l:four} Let $G$ be a pro-Lie group satisfying
hypothesis {\bf(*)} of {\rm Lemma \ref{l:two}}. Then $G_t$ is complete.
\end{Lemma}

\msk

\nin
We note right away that Lemma \ref{l:four} 
will prove hypothesis (H) and therefore
complete the proof of Theorem \ref{th:main}.

\msk

\nin
{\bf Proof} of Lemma \ref{l:four}. 
We let $f\colon G\to G_t=G/G_0$ be the
quotient morphism and consider a Cauchy filter $\cal C$ on $G_t$.
We have to show that $\cal C$ converges. By Lemma \ref{l:two}, 
${\cal M}(G)$
is cofinal in ${\cal N}(G)$. For each $N\in{\cal M}(G)$ let
$N^*=\overline{f(N)}$ and let 
$p_{N^*}\colon G\to G_t{/N^*}$ be the 
quotient morphism.
Then the image $p_{N^*}({\cal C})$ is a 
Cauchy filter in the Lie group $G_t/N^*$ and thus has a limit $g_N$.
Then $(g_N)_{N\in{\cal M}(G)}\in\prod_{N\in{\cal M}(G)}G_t/N^*$
is an element of $\lim_{N\in{\cal M}(G)}G_t/N^*$; indeed $\cal C$
has to converge to a point in the completion of $G_t$.
Now let $F_N=(p_{N^*}\circ f)^{-1}({g}_N)$. 
Then $\{F_N:N\in{\cal M}(G)\}$
is a filter basis consisting of cosets
 modulo $N_*\defi\overline{G_0N}$ of $G$.
We claim that $N_*\in{\cal M}(G)$. Indeed we have 
$N_*\in{\cal N}(G)$. Now we let $M$ be an open subgroup 
of $N_*\supseteq G_0$.
Then $G_0\subseteq M$, so $G_0N\subseteq MN$ and $MN$ 
is open-closed in $N_*=\overline{G_0N}$. 
Thus $N_*=MN$. So $MN/M \cong N/(N\cap M)$ is discrete
and then $M\cap N$ is open in $N$. But $N\in{\cal M}(G)$ then implies 
that $N_*/M\cong N/(M\cap N)$ is finite. This shows that 
$N_*\in{\cal M}(G)$
as claimed. Since $G$ is Z-compact by Lemma { \ref{l:three}}, 
we find an element 
$g\in\bigcap_{N\in{\cal M}(G)}F_N$. But then 
$p_{N^*}(f(g))={g}_N$ for
all $N\in{\cal M}(G)$ which implies that $f(g)=\lim{\cal C}$.
Thus every Cauchy filter in $G_t$ converges showing 
that $G_t$ is complete.
\qed

\msk 

\nin
This completes the proof of Theorem \ref{th:main}.

\bsk

\quad An inspection of \cite{probook} shows that
the following questions appear to be unsettled:

\msk

\nin
{\bf Question 1.}\quad For which pro-Lie groups $G$ is 
${\cal M}(G)$ cofinal in ${\cal N}(G)$?

\msk

For each of these groups $G$ we would know that $G$ is (isomorphic to)
the strict projective limit $\lim_{M\in{\cal M}(G)}G/M$.
In \cite{probook} this is proved of all almost
connected pro-Lie groups.

\ssk

\nin
Test examples are the nondiscrete pro-discrete groups 
$\Z^{(\N)}$ and $\Z^\N$
(see e.g. \cite{probook}, Example 4.4ff., Proposition 5.2).

\msk

\nin
{\bf Question 2.}\quad For which pro-Lie groups is the Z-topology 
compact?

\msk

In \cite{probook}, Proposition 4.27, this is shown for almost connected
pro-Lie groups, and here we have proved it for
those pro-Lie groups $G$ all of whose Lie group quotients are known to be 
almost connected.

\msk

Theorem \ref{th:main} suggests the following rather general question: 

\msk  

\nin {\bf Question 3.}\quad When is the projective limit of a projective system
of almost connected topological groups almost connected?

\msk
Theorem \ref{th:main} says that within the category of pro-Lie groups 
we have an affirmative answer
 for the projective system of all Lie group quotients.
See also some background information in \cite{probook} in and around
Theorem 1.27, p.~88.
\bsk

\nin
One remark is in order in the context of the Z-topology discussed in
\cite{probook} on pp.~197--203: 

In Exercise E4.2(i), {p.~202,} it is pointed out that
{ ${\cal M}(\Z)=\big\{\{0\},\Z\big\}$ and that therefore} the 
Z-topology on $\Z$, being the cofinite topology, is compact. {
Whereas the topology generated by the set of cosets $z+N$, 
$N\in{\cal N}(\Z)$, the set}  
of all subgroups of $\Z$,
  fails to be compact.

\bsk

Theorem \ref{th:main} will play a significant role 
in the authors' study \cite{dah} of weakly complete real or complex topological
algebras with identity, which will explore in detail  their relation 
to pro-Lie theory and  aims  for a systematic treatment of weakly
complete group algebras of topological groups and their representation
and duality theories.

\msk

{
\nin {\bf Acknowledgment.}\quad The authors thank the referee for his
swift, yet thorough contributions to the final form of this note.}

\end{document}